\documentclass[12pt]{amsart}
\usepackage{bbding}
\usepackage[margin=2.7cm]{geometry}

\usepackage{amsfonts}
\usepackage{mathrsfs}
\usepackage{cite}
\usepackage{tikz}

\newtheorem{theorem}{Theorem}[section]
\newtheorem{lemma}[theorem]{Lemma}

\newcommand{\Conv}[1]{\text{Conv }#1}
\newcommand{\Int}[1]{\text{Int }#1}

\title{The $L_{p}$ Minkowski problem for polytopes for negative $p$}
\author{guangxian zhu}
\address{Department of Mathematics, Polytechnic School of  Engineering, New York University, Brooklyn, NY 11201.}
\email{gz342@nyu.edu}

\begin{document}
\maketitle

\begin{abstract}
Existence of solutions to the $L_{p}$ Minkowski problem is proved for all $p<0$. For the critical case of $p=-n$, which is known as the centro-affine Minkowski problem, this paper contains the main result in \cite{ZHU} as a special case.
\end{abstract}

\section{Introduction}

A \emph{convex body}  in $n$-dimensional Euclidean space,
$\mathbb{R}^{n}$,  is a compact convex set that has
non-empty interior. If $p\in\mathbb{R}$ and $K$ is a convex body in $\mathbb{R}^{n}$ that
contains the origin in its interior, then the $L_{p}$ surface area
measure,
$S_{p}(K,\cdot)$, of $K$ is a Borel measure on the unit sphere, $S^{n-1}$, defined for
each Borel $\omega\subset S^{n-1}$ by
$$
S_{p}(K,\omega)=\int_{x\in\nu_{K}^{-1}(\omega)}(x\cdot\nu_{K}(x))^{1-p}d\mathcal{H}^{n-1}(x),
$$
where $\nu_{K}:\partial'K\rightarrow S^{n-1}$ is the Gauss map of $K$,
defined on $\partial'K$, the set of boundary points of $K$ that have a
unique outer unit normal, and $\mathcal{H}^{n-1}$ is
($n-1$)-dimensional Hausdorff measure.

The $L_{p}$ surface area measure was introduced by Lutwak \cite{LUT}. The $L_{p}$ surface area measure contains three important measures as special cases: the $L_{1}$ surface area measure is the classic surface area measure; the $L_{0}$ surface area measure is the cone-volume measure; the $L_{-n}$ surface area measure is the centro-affine surface area measure. Today, the $L_{p}$ surface area measure is a central notation in convex geometry analysis, and appeared in, e.g., \cite{BG, CG, H1, HLYZ1, HP1, HP2, HS1, HS2, HS3, HENK, LU1, ML3, LU2, LR, LUT, LO, LYZ1, LYZ2, LYZ3, LYZ4, LYZ5, LYZ6, LZ, MIN, NAO, NR, PW, RZ, SCH, FS2, FS4, FS5, WA, WE, YE}.

The following $L_{p}$ Minkowski problem that posed by Lutwak \cite{LUT} is considered as one of the most important problems in modern convex geometry analysis.\\ \\
\textbf{$L_{p}$ Minkowski problem:} \textit{Find necessary and sufficient conditions on a finite Borel measure $\mu$ on $S^{n-1}$ so that $\mu$
is the $L_{p}$ surface area measure of a convex body in $\mathbb{R}^{n}$}.\\

The associated partial differential equation for the $L_{p}$ Minkowski problem is the following Mong-Amp\`{e}re type equation: For a given positive function $f$ on the unit sphere, solve
\begin{equation}\tag{1.1}\label{Equation 1.1}
h^{1-p}\det(h_{ij}+h\delta_{ij})=f,
\end{equation}
where $h_{ij}$ is the covariant derivative of $h$ with respect to an orthonormal frame on $S^{n-1}$ and $\delta_{ij}$ is the Kronecker delta.

The solutions of the $L_{p}$ Minkowski problem have important applications to affine isoperimetric inequalities, see, e.g., Zhang \cite{Z}, Lutwak, Yang and Zhang \cite{LYZ4}, Ciachi, Lutwak, Yang and Zhang \cite{CLYZ}, Haberl and Schuster \cite{HS1, HS2, HS3}. The solutions to the $L_{p}$ Minkowski problem are also related with some important flows (see, e.g., \cite{AN1, AN2, ST1, ST2}).

When $p=1$, the $L_{p}$ Minkowski problem is the classical Minkowski problem. The existence and uniqueness for the solution of this problem was solved by Minkowski, Aleksandrov, and Fenchel and Jessen (see Schneider \cite{SCH} for references). Regularity of the Minkowski problem was studied by e.g., Caffarelli \cite{LC},  Cheng and Yau \cite{CY}, Nirenberg \cite{NIR} and Pogorelov \cite{POG}.

For $p\neq1$, the $L_{p}$ Minkowski problem was studied by, e.g., Lutwak \cite{LUT}, Lutwak and Oliker \cite{LO}, Lutwak, Yang and Zhang \cite{LYZ5}, Chou and Wang \cite{CW}, Guan and Lin \cite{GL}, Hug, Lutwak, Yang and Zhang \cite{HLYZ1}, B\"{o}r\"{o}czky, Heged\H{u}s and Zhu \cite{BHZ}, B\"{o}r\"{o}czky, Lutwak, Yang and Zhang \cite{BLYZ2, BLYZ}, Chen \cite{WC}, Dou and Zhu \cite{DZ}, Haberl, Lutwak, Yang and Zhang \cite{HLYZ1}, Huang, Liu and Xu \cite{HLX}, Jian, Lu and Wang \cite{JLW}, Jian and Wang \cite{JW}, Jiang, Wang and Wei \cite{JWW}, Lu and Wang \cite{LWA}, Stancu \cite{ST1, ST2}, Sun and Long \cite{SL} and Zhu \cite{ZHU2, ZHU, ZHU1}. Analogues of the Minkowski problems were studied in, e.g., \cite{COL, GAGE, GH, GG, GM, HMS, XIA}.

The uniqueness of solutions to the $L_{p}$ Minkowski for $p>1$ can be shown by applying the $L_{p}$ Minkowski inequality established by Lutwak \cite{LUT}. However, little is know about the $L_{p}$ Minkowski inequality for the case where $p<1$. This is one of the main reasons that most of the previous work on the $L_{p}$ Minkowski problem was limited to the case where $p>1.$

The critical case where $p=-n$ of the $L_{p}$ Minkowski problem is called the centro-affine Minkowski problem, which describes the centro-affine surface area measure. This problem is especially important due to the affine invariant of the partial differential equation (\ref{Equation 1.1}). It is known that the centro-affine Minkowski problem has connections with several important geometric problems (see, e.g., Jian and Wang \cite{JW} for reference). The centro-affine Minkowski problem was explicitly posed by Chow and Wang \cite{CW}. Recently, the centro-affine Minkowski problem was studied by Lu and Wang \cite{LWA} for rotationally symmetric case and was studied by Zhu \cite{ZHU} for discrete measures.

When $p<-n$, very few results are known for the $L_{p}$ Minkowski problem. So far as the author knows, in $\mathbb{R}^{2}$, the $L_{p}$ Minkowski problem for all $p<0$ was studied by Dou and Zhu \cite{DZ}, Sun and Long \cite{SL}. It is the aim of this paper to study the $L_{p}$ Minkowski problem for all $p<0$ and $n\geq2$.

It is know that the Minkowski problem and the $L_{p}$ Minkowski problem (for $p>1$) for arbitrary measures can be solved by an approximation argument by first solving the polytopal case (see, e.g., \cite{HLYZ2} or \cite{SCH} pp. 392-393). This is one of the reasons why the Minkowski problem and the $L_{p}$ Minkowski problem for polytopes are of great importance.

A \emph{polytope} in $\mathbb{R}^{n}$ is the convex hull of a finite set of points in $\mathbb{R}^{n}$ provided that it has positive $n$-dimensional volume. The convex hull of a subset of these points is called a \emph{facet} of the polytope if it lies entirely on the boundary of the polytope and has
positive ($n-1$)-dimensional volume. Let $P$ be a polytope which contains the origin in its interior with $N$ facets whose outer unit normals are
$u_{1},...,u_{N}$, and such that the facet with outer unit normal $u_{k}$ has area $a_{k}$ and distance $h_{k}$ from the origin for
all $k\in\{1,...,N\}$. Then,
$$
S_{p}(P,\cdot)=\sum_{k=1}^{N}h_{k}^{1-p}a_{k}\delta_{u_{k}}(\cdot).
$$
where $\delta_{u_{k}}$ denotes the delta measure that is concentrated at the point $u_{k}$.

A finite subset $U$ of $S^{n-1}$ is said to be \emph{in general position} if any $k$ elements of $U$, $1\le k\le n$,  are linearly independent.

In \cite{ZHU}, the author solved the centro-affine Minkowski problem for polytopes whose outer unit normals are in general position:\\ \\
\textbf{Theorem A.}\emph{ Let $\mu$ be a discrete measure on the unit sphere $S^{n-1}$. Then $\mu$ is the centro-affine surface area measure of a polytope whose outer unit normals are in general position if and only if the support of $\mu$ is in general position and not concentrated on a closed hemisphere.}\\

A linear subspace $X$ ($0<\dim X<n$) of  $\mathbb{R}^{n}$ is said to be \textit{essential} with respect to a Borel measure $\mu$ on $S^{n-1}$ if $X\cap{\rm supp}(\mu)$ is not concentrated on any closed hemisphere of $X\cap S^{n-1}$.

Obviously, if the support of a discrete measure $\mu$ is in general position, then the set of essential subspaces of $\mu$ is empty. On the other hand, in $\mathbb{R}^{n}$ ($n\geq3$), one can easily construct a discrete measure $\mu$ such that $\mu$ does not has essential subspace but the support of $\mu$ is not in general position. Therefore, the set of discrete measures whose supports are in general position is a subset of the set of discrete measures that do not have essential subspaces.

It is the aim of this paper to solve the $L_{p}$ Minkowski problem for discrete measures that do not have essential subspaces. Obviously, the following main theorem of this paper contains Theorem A as a special case.
\begin{theorem}
Let $p<0$ and $\mu$ be a discrete measure on the unit sphere $S^{n-1}$. Then $\mu$ is the $L_{p}$ surface area measure of a polytope whose $L_{p}$ surface area measure does not have essential subspace if and only if $\mu$ does not have essential subspace and not concentrated on a closed hemisphere.
\end{theorem}

\section{Preliminaries}

In this section, we standardize some notations and list some basic facts about convex bodies. For general references regarding convex bodies, see, e.g., \cite{GAR1, SCH, THO}.

The sets in this paper are subsets of the $n$-dimensional Euclidean space $\mathbb{R}^{n}$. For $x, y\in\mathbb{R}^{n}$, we write $x\cdot y$ for the standard inner product of $x$ and $y$, $|x|$ for the Euclidean norm of $x$, and $S^{n-1}$ for the unit sphere of $\mathbb{R}^{n}$.

Suppose $S$ is a subset of $\mathbb{R}^{n}$, then the positive hull, pos$(S)$, of $S$ is the set of all positive combinations of any finitely many elements of $S$. Let lin$(S)$ be the smallest linear subspace of $\mathbb{R}^{n}$ containing $S$. The diameter of a subset, $S$, of $\mathbb{R}^{n}$ is defined by
$$
d(S)=\max\{|x-y|: x,y\in S\}.
$$
The convex hull of a subset, $S$, of $\mathbb{R}^{n}$ is defined by
$$
\Conv(S)=\{\lambda x+(1-\lambda)y: 0\leq\lambda\leq1\textrm{ and }x, y\in S\}.
$$

For convex bodies $K_{1}, K_{2}$ in $\mathbb{R}^{n}$ and $s_{1}, s_{2}\geq0$, the Minkowski combination is defined by
$$
s_{1}K_{1}+s_{2}K_{2}=\{s_{1}x_{1}+s_{2}x_{2}: x_{1}\in K_{1}, x_{2}\in K_{2}\}.
$$
The \emph{support function}
$h_{K}: \mathbb{R}^{n}\rightarrow\mathbb{R}$ of a convex body
$K$ is defined, for $x\in\mathbb{R}^{n}$, by
$$
h(K,x)=\max\{x\cdot y: y\in K\}.
$$
Obviously, for $s\geq0$ and $x\in\mathbb{R}^{n}$,
$$
h(sK,x)=h(K,sx)=sh(K,x).
$$

If $K$ is a convex body in $\mathbb{R}^{n}$ and $u\in S^{n-1}$, then
the \emph{support set} $F(K,u)$ of $K$ in direction $u$ is defined by
$$
F(K,u)=K\cap\{x\in\mathbb{R}^{n}:x\cdot u=h(K,u)\}.
$$

The \emph{Hausdorff distance} of two convex bodies $K_{1}, K_{2}$ in
$\mathbb{R}^{n}$ is defined by
$$
\delta(K_{1},K_{2})=\inf\{t\geq0: K_{1}\subset K_{2}+tB^{n}, K_{2}\subset K_{1}+tB^{n}\},
$$
where $B^{n}$ is the unit ball.

Let $\mathcal{P}$ be the set of polytopes in $\mathbb{R}^{n}$. If the unit vectors $u_{1},...,u_{N}$ are not concentrated on a closed hemisphere, let $\mathcal{P}(u_{1},...,u_{N})$ be the subset of $\mathcal{P}$ such that a polytope $P\in\mathcal{P}(u_{1},...,u_{N})$ if the the set of the outer unit normals of $P$ is a subset of $\{u_{1},...,u_{N}\}$. Let $\mathcal{P}_{N}(u_{1},...,u_{N})$ be the subset of $\mathcal{P}(u_{1},...,u_{N})$ such that a
polytope $P\in\mathcal{P}_{N}(u_{1},...,u_{N})$ if, $P\in\mathcal{P}(u_{1},...,u_{N})$, and $P$ has exactly $N$ facets.

\section{An extremal problem related to the $L_{p}$ Minkowski problem}

Suppose $p<0$, $\alpha_{1},...,\alpha_{N}>0$, the unit vectors $u_{1},...,u_{N}$ are not concentrated on a closed hemisphere, and $P\in\mathcal{P}(u_{1},...,u_{N})$. Define the function, $\Phi_{P}: \Int(P)\rightarrow\mathbb{R}$, by
$$
\Phi_{P}(\xi)=\sum_{k=1}^{N}\alpha_{k}(h(P,u_{k})-\xi\cdot u_{k})^{p}.
$$

In this section, we study the extremal problem
\begin{equation}\tag{3.1}\label{Equation 3.1}
\sup\{\inf_{\xi\in\Int(Q)}\Phi_{Q}(\xi): Q\in\mathcal{P}(u_{1},...,u_{N})\textmd{ and }V(Q)=1\}.
\end{equation}

The main purpose of this section is to prove that a dilation of the solution to problem (3.1) solves the corresponding $L_{p}$ Minkowski problem.

\begin{lemma}\label{Lemma 3.1}
 If $p<0$, $\alpha_{1},...,\alpha_{N}>0$, the unit vectors $u_{1},...,u_{N}$ are not concentrated on a closed hemisphere and $P\in\mathcal{P}(u_{1},...,u_{N})$, then there exists a unique $\xi(P)\in\Int(P)$ such that
 $$
 \Phi_{P}(\xi(P))=\inf_{\xi\in\Int(P)}\Phi_{P}(\xi).
 $$
\end{lemma}
\begin{proof}
Since $p<0$, the function $f(t)=t^{p}$ is strictly convex on $(0,+\infty)$. Hence, for $0<\lambda<1$
and $\xi_{1}, \xi_{2}\in\Int(P)$,
\begin{equation*}
\begin{split}
\lambda\Phi_{P}(\xi_{1})+(1-\lambda)\Phi_{P}(\xi_{2})&=\lambda\sum_{k=1}^{N}\alpha_{k}(h(P,u_{k})-\xi_{1}\cdot
u_{k})^{p}+(1-\lambda)\sum_{k=1}^{N}\alpha_{k}(h(P,u_{k})-\xi_{2}\cdot
u_{k})^{p}\\
&=\sum_{k=1}^{N}\alpha_{k}\left[\lambda(h(P,u_{k})-\xi_{1}\cdot
u_{k})^{p}+(1-\lambda)(h(P,u_{k})-\xi_{2}\cdot u_{k})^{p}\right]\\
&\geq\sum_{k=1}^{N}\alpha_{k}\left[h(P,u_{k})-(\lambda\xi_{1}+(1-\lambda)\xi_{2})\cdot
u_{k}\right]^{p}\\
&=\Phi_{P}(\lambda\xi_{1}+(1-\lambda)\xi_{2}).
\end{split}
\end{equation*}
Equality hold if and only if $\xi_{1}\cdot u_{k}=\xi_{2}\cdot u_{k}$
for all $k=1,...,N$. Since $u_{1},...,u_{N}$ are not concentrated on a closed hemisphere, $\mathbb{R}^{n}=\textmd{lin}\{u_{1},...,u_{N}\}$. Thus, $\xi_{1}=\xi_{2}$. Hence, $\Phi_{P}$ is strictly convex on $\Int(P)$.

From the fact that $P\in\mathcal{P}(u_{1},...,u_{N})$, we have, for any $x\in\partial P$,
there exists a $u_{i_{0}}\in\{u_{1},...,u_{N}\}$ such that
$$
h(P,u_{i_{0}})=x\cdot u_{i_{0}}.
$$
Thus, $\Phi_{P}(\xi)\rightarrow\infty$ whenever $\xi\in\Int(P)$ and $\xi\rightarrow x$. Therefore, there exists a unique interior point
$\xi(P)$ of $P$
such that
$$
\Phi_{P}(\xi(P))=\inf_{\xi\in\Int(P)}\Phi_{P}(\xi).
$$
\end{proof}

Obviously, for $\lambda>0$ and $P\in\mathcal{P}(u_{1},...,u_{N})$,
\begin{equation}\tag{3.2}\label{Equation 3.2}
\xi(\lambda P)=\lambda\xi(P),
\end{equation}
and if $P_{i}\in\mathcal{P}(u_{1},...,u_{N})$ and $P_{i}$
converges to a polytope $P$, then
$P\in\mathcal{P}(u_{1},...,u_{N})$.

\begin{lemma}
If $p<0$, $\alpha_{1},...,\alpha_{N}>0$, the unit vectors $u_{1},...,u_{N}$ are not contained in a closed hemisphere, $P_{i}\in\mathcal{P}(u_{1},...,u_{N})$, and $P_{i}$ converges to a polytope $P$, then
$\lim_{i\rightarrow\infty}\xi(P_{i})=\xi(P)$ and
$$
\lim_{i\rightarrow\infty}\Phi_{P_{i}}(\xi(P_{i}))=\Phi_{P}(\xi(P)).
$$
\end{lemma}
\begin{proof}
Since $P_{i}$ converges to $P$ and $\xi(P_{i})\in\Int(P_{i})$, $\xi(P_{i})$ is bounded. Let $\xi_{0}$ be the limit point of a subsequence, $\xi(P_{i_{j}})$, of $\xi(P_{i})$. We claim that $\xi_{0}\in\Int(P)$. Otherwise, $\xi_{0}$ is a boundary point of $P$ with $\lim_{j\rightarrow\infty}\Phi_{P_{i_{j}}}(\xi_{P_{i_{j}}})=\infty$, which contradicts the fact that
\begin{equation}\tag{3.3}\label{Equation 3.3}
\overline{\lim_{j\rightarrow\infty}}\Phi_{P_{i_{j}}}(\xi(P_{i_{j}}))\leq\overline{\lim_{j\rightarrow\infty}}\Phi_{P_{i_{j}}}(\xi(P))=\Phi(\xi(P))<\infty.
\end{equation}

We claim that $\xi_{0}=\xi(P)$. Otherwise,
\begin{equation*}
\begin{split}
\lim_{j\rightarrow\infty}\Phi_{P_{i_{j}}}(\xi(P_{i_{j}}))&=\Phi_{P}(\xi_{0})\\
&>\Phi_{P}(\xi(P))\\
&=\lim_{j\rightarrow\infty}\Phi_{P_{i_{j}}}(\xi(P)).
\end{split}
\end{equation*}
This contradicts the fact that
$$
\Phi_{P_{i_{j}}}(\xi(P_{i_{j}}))\leq\Phi_{P_{i_{j}}}(\xi(P)).
$$
Hence, $\lim_{i\rightarrow\infty}\xi(P_{i})=\xi(P)$ and
$$
\lim_{i\rightarrow\infty}\Phi_{P_{i}}(\xi(P_{i}))=\Phi_{P}(\xi(P)).
$$
\end{proof}

\begin{lemma}
If $p<0$, $\alpha_{1},...,\alpha_{N}>0$, the unit vectors $u_{1},...,u_{N}$ are not concentrated on a closed hemisphere and $P\in\mathcal{P}(u_{1},...,u_{N})$, then
$$
\sum_{k=1}^{N}\alpha_{k}\frac{u_{k}}{[h(P,u_{k})-\xi(P)\cdot u_{k}]^{1-p}}=0.
$$
\end{lemma}
\begin{proof}
Define $f:$Int$(P)\rightarrow\mathbb{R}^{n}$ by
$$
f(x)=\sum_{k=1}^{N}\alpha_{k}(h(P,u_{k})-x\cdot u_{k})^{p}.
$$

By conditions,
$$
f(\xi(P))=\inf_{x\in Int(P)}f(x).
$$
Thus,
$$
\sum_{k=1}^{N}\alpha_{k}\frac{u_{k,i}}{[h(P,u_{k})-\xi(P)\cdot u_{k}]^{1-p}}=0,
$$
for all $i=1,...,n$, where $u_{k}=(u_{k,1},...,u_{k,n})^{T}$. Therefore,
$$
\sum_{k=1}^{N}\alpha_{k}\frac{u_{k}}{[h(P,u_{k})-\xi(P)\cdot u_{k}]^{1-p}}=0.
$$
\end{proof}
\begin{lemma}
 Suppose $p<0$, $\alpha_{1},..,\alpha_{N}>0$, the unit vectors $u_{1},...,u_{N}$ are not concentrated on a closed hemisphere, and there exists a $P\in\mathcal{P}_{N}(u_{1},...,u_{N})$ with
$\xi(P)=o$, $V(P)=1$ such that
$$
\Phi_{P}(o)=\sup\left\{\inf_{\xi\in\Int (Q)}\Phi_{Q}(\xi): Q\in\mathcal{P}(u_{1},...,u_{N})\textmd{ and }V(Q)=1\right\}.
$$
Then,
$$
S_{p}(P_{0},\cdot)=\sum_{k=1}^{N}\alpha_{k}\delta_{u_{k}}(\cdot),
$$
where $P_{0}=\left(\sum_{j=1}^{N}\alpha_{j}h(P,u_{j})^{p}/n\right)^{\frac{1}{n-p}}P.$
\end{lemma}
\begin{proof}
By conditions, there exists a polytope $P\in\mathcal{P}_{N}(u_{1},...,u_{N})$ with $\xi(P)=o$ and $V(P)=1$ such that
$$
\Phi_{P}(o)=\sup\{\inf_{\xi\in\Int(Q)}\Phi_{Q}(\xi): Q\in\mathcal{P}(u_{1},...,u_{N})\textmd{ and }V(Q)=1\},
$$
where $\Phi_{Q}(\xi)=\sum_{k=1}^{N}\alpha_{k}(h(Q,u_{k})-\xi\cdot u_{k})^{p}$.

For $\tau_{1},...,\tau_{N}\in\mathbb{R}$, choose $|t|$ small enough so that the polytope $P_{t}$ defined by
$$
P_{t}=\bigcap_{i=1}^{N}\left\{x: x\cdot u_{i}\leq h(P,u_{i})+t\tau_{i}\right\}
$$
has exactly $N$ facets. By \cite{SCH} (Lemma 7.5.3),
$$
\frac{\partial V(P_{t})}{\partial t}=\sum_{i=1}^{N}\tau_{i}a_{i},
$$
where $a_{i}$ is the area of $F(P,u_{i})$. Let $\lambda(t)=V(P_{t})^{-\frac{1}{n}}$, then $\lambda(t)P_{t}\in\mathcal{P}_{N}^{n}(u_{1},...,u_{N})$, $V(\lambda(t)P_{t})=1$ and
\begin{equation}\tag{3.4}\label{Equation 3.4}
\lambda'(0)=-\frac{1}{n}\sum_{i=1}^{N}\tau_{i}S_{i}.
\end{equation}

Define $\xi(t):=\xi(\lambda(t)P_{t})$, and
\begin{equation}\tag{3.5}\label{Equation 3.5}
\begin{split}
\Phi(t)&:=\min_{\xi\in\lambda(t)P_{t}}\sum_{k=1}^{N}\alpha_{k}(\lambda(t)h(P_{t},u_{k})-\xi\cdot u_{k})^{p}\\
&=\sum_{k=1}^{N}\alpha_{k}(\lambda(t)h(P_{t},u_{k})-\xi(t)\cdot u_{k})^{p}.
\end{split}
\end{equation}

It follows from Lemma 3.3 that
\begin{equation*}
\sum_{k=1}^{N}\alpha_{k}\frac{u_{k,i}}{[\lambda(t)h(P_{t},u_{k})-\xi(t)\cdot u_{k}]^{1-p}}=0,
\end{equation*}
for $i=1,...,n,$ where $u_{k}=(u_{k,1},...,u_{k,n})^{T}$. In addition, since $\xi(P)$ is the origin,
\begin{equation}\tag{3.6}\label{Equation 3.6}
\sum_{k=1}^{N}\alpha_{k}\frac{u_{k}}{h(P,u_{k})^{1-p}}=0.
\end{equation}

Let $F=(F_{1},...,F_{n})$ be a function from an open neighbourhood of the origin in $\mathbb{R}^{n+1}$ to $\mathbb{R}^{n}$ such that
$$
F_{i}(t,\xi_{1},...,\xi_{n})=\sum_{k=1}^{N}\alpha_{k}\frac{u_{k,i}}{[\lambda(t)h(P_{t},u_{k})-(\xi_{1}u_{k,1}+...+\xi_{n}u_{k,n})]^{1-p}}
$$
for $i=1,...,n.$ Then,
$$
\frac{\partial F_{i}}{\partial t}\bigg|_{(t,\xi_{1},...,\xi_{n})}=\sum_{k=1}^{N}\frac{(p-1)\alpha_{k}u_{k,i}\left[\lambda'(t)h(P_{t},u_{k})+\lambda(t)\tau_{k}\right]}{\left[\lambda(t)h(P_{t}, u_{k})-(\xi_{1}u_{k,1}+...+\xi_{n}u_{k,n})\right]^{2-p}},
$$
$$
\frac{\partial F_{i}}{\partial\xi_{j}}\bigg|_{(t,\xi_{1},...,\xi_{n})}=\sum_{k=1}^{N}\frac{(1-p)\alpha_{k}u_{k,i}u_{k,j}}{\left[\lambda(t)h(P_{t}, u_{k})-(\xi_{1}u_{k,1}+...+\xi_{n}u_{k,n})\right]^{2-p}}
$$
are continuous on a small neighbourhood of $(0,0,...,0)$ with
$$
\left(\frac{\partial F}{\partial\xi}\bigg|_{(0,...,0)}\right)_{n\times n}=\sum_{k=1}^{N}\frac{(1-p)\alpha_{k}}{h(P,u_{k})^{2-p}}u_{k}\cdot
u_{k}^{T},
$$
where $u_{k}u_{k}^{T}$ is an $n\times n$ matrix.

Since $u_{1},...,u_{N}$ are not contained in a closed hemisphere, $\mathbb{R}^{n}=\textmd{lin}\{u_{1},...,u_{N}\}$. Thus, for any $x\in\mathbb{R}^{n}$ with $x\neq0$, there exists a $u_{i_{0}}\in\{u_{1},...,u_{N}\}$ such that
$u_{i_{0}}\cdot x\neq0$. Then,
\begin{equation*}
\begin{split}
x^{T}\cdot\left(\sum_{k=1}^{N}\frac{(1-p)\alpha_{k}}{h(P,u_{k})^{2-p}}u_{k}\cdot u_{k}^{T}\right)\cdot x&=\sum_{k=1}^{N}\frac{(1-p)\alpha_{k}}{h(P,u_{k})^{2-p}}(x\cdot u_{k})^{2}\\
&\geq\frac{(1-p)\alpha_{i_{0}}}{h(P,u_{i_{0}})^{2-p}}(x\cdot u_{i_{0}})^{2}>0.
\end{split}
\end{equation*}
Therefore, $(\frac{\partial F}{\partial\xi}\big|_{(0,...,0)})$ is positive defined. By this, the fact that $F_{i}(0,...,0)=0$ for all $i=1,...,n$, the fact that $\frac{\partial F_{i}}{\partial\xi_{j}}$ is continuous on a neighbourhood of $(0,0,...,0)$ for all $0\leq i, j\leq n$ and the implicit function theorem, we have
$$
\xi'(0)=(\xi_{1}'(0),...,\xi_{n}'(0))
$$
exists.

From the fact that $\Phi(0)$ is an extreme value of $\Phi(t)$ (in Equation (\ref{Equation 3.5})), Equation (\ref{Equation 3.4}) and Equation (\ref{Equation 3.6}), we have
\begin{equation*}
\begin{split}
0&=\Phi'(0)/p\\
&=\sum_{k=1}^{N}\alpha_{k}h(P,u_{k})^{p-1}\left(\lambda'(0)h(P,u_{k})+\tau_{k}-\xi'(0)\cdot u_{k}\right)\\
&=\sum_{k=1}^{N}\alpha_{k}h(P,u_{k})^{p-1}\left[-\frac{1}{n}\left(\sum_{i=1}^{N}a_{i}\tau_{i}\right)h(P,u_{k})+\tau_{k}\right]-\xi'(0)\cdot\left[\sum_{k=1}^{N}
\alpha_{k}\frac{u_{k}}{h(P,u_{k})^{1-p}}\right]\\
&=\sum_{k=1}^{N}\alpha_{k}h(P,u_{k})^{p-1}\tau_{k}-\left(\sum_{i=1}^{N}a_{i}\tau_{i}\right)\frac{\sum_{k=1}^{N}\alpha_{k}h(P,u_{k})^{p}}{n}\\
&=\sum_{k=1}^{N}\left(\alpha_{k}h(P,u_{k})^{p-1}-\frac{\sum_{j=1}^{N}\alpha_{j}h(P,u_{j})^{p}}{n}a_{k}\right)\tau_{k}.
\end{split}
\end{equation*}
Since $\tau_{1},...,\tau_{N}$ are arbitrary,
$$
\frac{\sum_{j=1}^{N}\alpha_{j}h(P,u_{j})^{p}}{n}h(P,u_{k})^{1-p}a_{k}=\alpha_{k},
$$
for all $k=1,...,N$. By letting
$$
P_{0}=\left(\frac{\sum_{j=1}^{N}\alpha_{j}h(P,u_{j})^{p}}{n}\right)^{\frac{1}{n-p}}P,
$$
we have
$$
S_{p}(P_{0},\cdot)=\sum_{k=1}^{N}\alpha_{k}\delta_{u_{k}}(\cdot).
$$
\end{proof}

\section{The proof of the main theorem}

In this section, we prove the main theorem of this paper.

The following lemmas will be needed.

\begin{lemma}
Let $\{h_{1j}\}_{j=1}^{\infty},...,\{h_{Nj}\}_{j=1}^{\infty}$ be $N$ $(N\geq2)$ sequences of real numbers. Then, there exists a subsequence, $\{j_{n}\}_{n=1}^{\infty}$, of $\mathbb{N}$ and a rearrangement, $i_{1},...,i_{N}$, of $1,...,N$ such that
$$
h_{i_{1}j_{n}}\leq h_{i_{2}j_{n}}\leq...\leq h_{i_{N}j_{n}},
$$
for all $n\in\mathbb{N}$.
\end{lemma}
\begin{proof}
For each fixed $j$, the number of the possible order (from small to big) of $h_{1j},...,h_{Nj}$ is $N!$. Therefore, there exists a subsequence, $\{j_{n}\}_{n=1}^{\infty}$, of $\mathbb{N}$ and a rearrangement, $i_{1},...,i_{N}$, of $1,...,N$ such that
$$
h_{i_{1}j_{n}}\leq h_{i_{2}j_{n}}\leq...\leq h_{i_{N}j_{n}},
$$
for all $n\in\mathbb{N}$.
\end{proof}
\begin{lemma}
Suppose the unit vectors $u_{1},...,u_{N}$ are not concentrated on a closed hemisphere, and for any subspace, $X$, of $\mathbb{R}^{n}$ with $1\leq\dim X\leq n-1,$ $\{u_{1},...,u_{N}\}\cap X$ is concentrated on a closed hemisphere of $S^{n-1}\cap X$. If $P_{m}$ is a sequence of polytopes with $V(P_{m})=1$, $o\in\Int(P_{m})$ and $P_{m}\in\mathcal{P}(u_{1},...,u_{N})$, then $P_{m}$ is bounded.
\end{lemma}
\begin{proof}
We only need to prove that if the diameter, $d(P_{i})$, of $P_{i}$ is not bounded, then there exists a subspace, $X$, of $\mathbb{R}^{n}$ with $1\leq\dim(X)\leq n-1$ and $\{u_{1},...,u_{N}\}\cap X$ is not concentrated on a closed hemisphere of $S^{n-1}\cap X$.

Let $\mu$ be a discrete measure on the unit sphere such that supp$(\mu)=\{u_{1},...,u_{N}\}$, $\mu(u_{i})=\alpha_{i}>0$ for $1\leq i\leq N$. Obviously, we only need to prove the lemma under the condition that $\xi(P_{m})=o$ for all $m\in\mathbb{N}$.

By Lemma 4.1, we may assume that
\begin{equation}\tag{4.0}\label{Equation 4.0}
h(P_{m},u_{1})\leq...\leq h(P_{m},u_{N}).
\end{equation}
By this and the condition that $V(P_{m})=1$ and $\lim_{m\rightarrow\infty}d(P_{m})=\infty$,
$$
\lim_{m\rightarrow\infty}h(P_{m},u_{1})=0\textmd{ and }\lim_{m\rightarrow\infty}h(P_{m},u_{N})=\infty.
$$
By this and (\ref{Equation 4.0}), there exists an $i_{0}$ ($1\leq i_{0}\leq N$) such that
\begin{equation}\tag{4.1}\label{Equation 4.1}
\overline{\lim}_{m\rightarrow\infty}\frac{h(P_{m},u_{i_{0}})}{h(P_{m},u_{1})}=\infty,
\end{equation}
and for $1\leq i\leq i_{0}-1$
\begin{equation}\tag{4.2}\label{Equation 4.2}
\overline{\lim}_{m\rightarrow\infty}\frac{h(P_{m},u_{i})}{h(P_{m},u_{1})}
\end{equation}
exists and equals to a positive number.

Let
$$
\Sigma=\textmd{pos}\{u_{1},...,u_{i_{0}-1}\}
$$
and
$$
\Sigma^{*}=\{x\in\mathbb{R}^{n}: x\cdot u_{i}\leq0\textmd{ for all }1\leq i\leq i_{0}-1\}.
$$

Let $1\leq j\leq i_{0}-1$ and $x\in\Sigma^{*}\cap S^{n-1}$. From the condition that $\xi(P_{m})$ is the origin and Lemma 3.3, we have
$$
\sum_{i=0}^{N}\frac{\alpha_{i}(x\cdot u_{i})}{[h(P_{m},u_{i})]^{1-p}}=0.
$$
By this and the fact that $x\in\Sigma^{*}\cap S^{n-1}$,
\begin{equation*}
\begin{split}
0&\geq\alpha_{j}(x\cdot u_{j})\\
&=-\sum_{i\neq j}\left[\frac{h(P_{m},u_{j})}{h(P_{m},u_{i})}\right]^{1-p}\alpha_{i}(x\cdot u_{i})\\
&\geq\sum_{i\geq i_{0}}\left[\frac{h(P_{m},u_{j})}{h(P_{m},u_{i})}\right]^{1-p}\alpha_{i}(x\cdot u_{i})\\
&\geq-\sum_{i\geq i_{0}}\left[\frac{h(P_{m},u_{j})}{h(P_{m},u_{i})}\right]^{1-p}\alpha_{i}
\end{split}
\end{equation*}
By this, (\ref{Equation 4.0}), (\ref{Equation 4.1}) and (\ref{Equation 4.2}), $\alpha_{j}(x\cdot u_{j})$ is no bigger than 0 and no less than any negative number. Hence,
$$
x\cdot u_{j}=0
$$
for all $j=1,...,i_{0}-1$ and $x\in\Sigma^{*}\cap S^{n-1}$. Thus,
\begin{equation}\tag{4.3}\label{Equation 4.3}
\Sigma^{*}\cap\textmd{lin}\{u_{1},...,u_{i_{0}-1}\}=\{0\}.
\end{equation}

Obviously, $\{u_{1},...,u_{i_{0}-1}\}$ is not concentrated on a closed hemisphere of $S^{n-1}\cap\textmd{lin}\{u_{1},...,u_{i_{0}-1}\}$. Otherwise, there exists an $x_{0}\in\textmd{lin}\{u_{1},...,u_{i_{0}-1}\}$ with $x_{0}\neq0$ such that $x_{0}\cdot u_{i}\leq0$ for all $1\leq i\leq i_{0}-1$. This contradicts with (\ref{Equation 4.3}).

We next prove that
$$
\textmd{lin}\{u_{1},...,u_{i_{0}-1}\}\neq\mathbb{R}^{n}.
$$
Otherwise, from the fact that $u_{1},...,u_{i_{0}-1}$ are not concentrated on a closed hemisphere of
$$
\textmd{lin}\{u_{1},...,u_{i_{0}-1}\}\cap S^{n-1},
$$
we have, the convex hull of $\{u_{1},...,u_{i_{0}-1}\}$ (denoted by $Q$) is a polytope in $\mathbb{R}^{n}$ and contains the origin as an interior. Let $F$ be a facet of $Q$ such that $\{su_{i_{0}}: s>0\}\cap F\neq\emptyset$. Since $F$ is the union of finite $(n-1)-$dimensional simplexes and the vertexes of these simplexes are subsets of $\{u_{1},...,u_{i_{0}-1}\}$, there exists a subset, $\{u_{i_{1}},...,u_{i_{n}}\}$, of $\{u_{1},...,u_{i_{0}-1}\}$ such that
$$
u_{i_{0}}\in\textmd{pos}\{u_{i_{1}},...,u_{i_{n}}\}.
$$

Since $o\in\Int(Q)$, there exists $r>0$ such that $rB^{n}\subset Q$. Choose $t>0$ such that $tu\in F\cap\textmd{pos}\{u_{i_{1}},...,u_{i_{n}}\}$. Then,
$$
tu=\beta_{i_{1}}u_{i_{1}}+...+\beta_{i_{n}}u_{i_{n}},
$$
where $\beta_{i_{1}},...,\beta_{i_{n}}\geq0$ with $\beta_{i_{1}}+...+\beta_{i_{n}}=1.$ If we let $a_{i_{j}}=\beta_{i_{j}}/t$ for $j=1,...,n$, we have
$$
u=a_{i_{1}}u_{i_{1}}+...+a_{i_{n}}u_{i_{n}}.
$$
Obviously, $a_{i_{j}}\geq0$ with
$$
a_{i_{j}}=\beta_{i_{j}}/t\leq1/r
$$
for all $j=1,...,n.$
Hence,
\begin{equation*}
\begin{split}
h(P_{m}, u_{i_{0}})&=h(P_{m},a_{i_{1}}u_{i_{1}}+...+a_{i_{n}}u_{i_{n}})\\
&\leq a_{i_{1}}h(P_{m},u_{i_{1}})+...+a_{i_{n}}h(P_{m},u_{i_{n}})\\
&\leq\frac{1}{r}\left[h(P_{m},u_{i_{1}})+...+h(P_{m},u_{i_{n}})\right],
\end{split}
\end{equation*}
for all $m\in\mathbb{N}$. This contradicts (\ref{Equation 4.1}) and (\ref{Equation 4.2}). Therefore,
$$
\textmd{lin}\{u_{1},...,u_{i_{0}-1}\}\neq\mathbb{R}^{n}.
$$

Let $X=\textmd{lin}\{u_{1},...,u_{i_{0}-1}\}$. Then, $1\leq\dim X\leq n-1$ but $\{u_{1},...,u_{N}\}\cap X=\{u_{1},...,u_{i_{0}-1}\}$ is not concentrated on a closed hemisphere of $S^{n-1}\cap X$, which contradicts the conditions of this lemma. Therefore, $d(P_{m})$ is bounded.
\end{proof}

The following lemmas will be needed (see, e.g., \cite{ZHU1}).
\begin{lemma}\label{Lemma 3.5}
If $P$ is a polytope in $\mathbb{R}^{n}$ and $v_{0}\in S^{n-1}$ with $V_{n-1}(F(P,v_{0}))=0$, then there exists a $\delta_{0}>0$ such that for $0\leq\delta<\delta_{0}$
$$
V(P\cap\{x: x\cdot v_{0}\geq h(P,v_{0})-\delta\})=c_{n}\delta^{n}+...+c_{2}\delta^{2},
$$
where $c_{n},...,c_{2}$ are constants that depend on $P$ and $v_{0}$.
\end{lemma}

\begin{lemma}
Suppose $p<0$, $\alpha_{1},...,\alpha_{N}>0$, and the unit vectors $u_{1},...,u_{N}$ are not concentrated on a hemisphere. If for any subspace $X$ with $1\leq\dim X\leq n-1$, $\{u_{1},...,u_{N}\}\cap X$ is always concentrated on a closed hemisphere of $S^{n-1}\cap X$,
then there exists a $P\in\mathcal{P}_{N}(u_{1},...,u_{N})$ such that $\xi(P)=o$, $V(P)=1$, and
$$
\Phi_{P}(o)=\sup\{\inf_{\xi\in\Int(Q)}\Phi_{Q}(\xi): Q\in\mathcal{P}(u_{1},...,u_{N})\textmd{ and }V(Q)=1\},
$$
where $\Phi_{Q}(\xi)=\sum_{k=1}^{N}\alpha_{k}(h(Q,u_{k})-\xi\cdot u_{k})^{p}$.
\end{lemma}
\begin{proof}
Obviously, for $P, Q\in\mathcal{P}(u_{1},...,u_{N})$, if there exists a $x\in\mathbb{R}^{n}$ such that $P=Q+x$, then
$$
\Phi_{P}(\xi(P))=\Phi_{Q}(\xi(Q)).
$$
Thus, we can choose a sequence of polytopes $P_{i}\in\mathcal{P}(u_{1},...,u_{N})$ with $\xi(P_{i})=o$ and $V(P_{i})=1$ such that $\Phi_{P_{i}}(o)$ converges
to
$$
\sup\{\inf_{\xi\in\Int(Q)}\Phi_{Q}(\xi): Q\in\mathcal{P}(u_{1},...,u_{N})\textmd{ and }V(Q)=1\}.
$$

By the conditions of this lemma and Lemma 4.2, $P_{i}$ is bounded. From the Blaschke selection theorem, there exists a subsequence of $P_{i}$ that converges to a
polytope $P$ such that $P\in\mathcal{P}(u_{1},...,u_{N})$, $V(P)=1$, $\xi(P)=o$ and
\begin{equation}\tag{4.4}\label{Equation 4.4}
\Phi_{P}(o)=\sup\{\inf_{\xi\in\Int(Q)}\Phi_{Q}(\xi): Q\in\mathcal{P}(u_{1},...,u_{N})\textmd{ and }V(Q)=1\}.
\end{equation}

We claim that $F(P,u_{i})$ are facets for all $i=1,...,N$. Otherwise, there exists an $i_{0}\in\{1,...,N\}$ such that
$$
F(P,u_{i_{0}})
$$
is not a facet of $P$.

Choose $\delta>0$ small enough so that the polytope
$$
P_{\delta}=P\cap\{x: x\cdot u_{i_{0}}\leq h(P,u_{i_{0}})-\delta\}\in\mathcal{P}(u_{1},...,u_{N}),
$$
and (by Lemma \ref{Lemma 3.5})
$$
V(P_{\delta})=1-(c_{n}\delta^{n}+...+c_{2}\delta^{2}),
$$
where $c_{n},...,c_{2}$ are constants that depend on $P$ and direction $u_{i_{0}}$.

From Lemma 3.2, for any $\delta_{i}\rightarrow0$ it always true that $\xi(P_{\delta_{i}})\rightarrow o$. We have,
$$
\lim_{\delta\rightarrow0}\xi(P_{\delta})=o.
$$
Let $\delta$ be small enough so that $h(P,u_{k})>\xi(P_{\delta})\cdot u_{k}+\delta$ for all $k\in\{1,...,N\}$, and let
$$
\lambda=V(P_{\delta})^{-\frac{1}{n}}=(1-(c_{n}\delta^{n}+...+c_{2}\delta^{2}))^{-\frac{1}{n}}.
$$
From this and Equation (\ref{Equation 3.2}), we have
\begin{equation}\tag{4.5}\label{Equation 4.5}
\begin{split}
\Phi_{\lambda P_{\delta}}(\xi(\lambda P_{\delta}))&=\sum_{k=1}^{N}\alpha_{k}\big(h(\lambda P_{\delta},u_{k})-\xi(\lambda
P_{\delta})\cdot
u_{k}\big)^{p}\\
&=\lambda^{p}\sum_{k=1}^{N}\alpha_{k}\big(h(P_{\delta},u_{k})-\xi(P_{\delta})\cdot u_{k}\big)^{p}\\
&=\lambda^{p}\sum_{k=1}^{N}\alpha_{k}\big(h(P,u_{k})-\xi(P_{\delta})\cdot
u_{k}\big)^{p}-\alpha_{i_{0}}\lambda^{p}\big(h(P,u_{i_{0}})-\xi(P_{\delta})\cdot
u_{i_{0}}\big)^{p}\\
&\qquad+\alpha_{i_{0}}\lambda^{p}\big(h(P,u_{i_{0}})-\xi(P_{\delta})\cdot u_{i_{0}}-\delta\big)^{p}\\
&=\sum_{k=1}^{N}\alpha_{k}\big(h(P,u_{k})-\xi(P_{\delta})\cdot
u_{k}\big)^{p}+(\lambda^{p}-1)\sum_{k=1}^{N}\alpha_{k}\big(h(P,u_{k})-\xi(P_{\delta})\cdot u_{k}\big)^{p}\\
&\qquad+\alpha_{i_{0}}\lambda^{p}\bigg[\big(h(P,u_{i_{0}})-\xi(P_{\delta})\cdot
u_{i_{0}}-\delta\big)^{p}-\big(h(P,u_{i_{0}})-\xi(P_{\delta})\cdot u_{i_{0}}\big)^{p}\bigg]\\
&=\Phi_{P}(\xi(P_{\delta}))+B(\delta),
\end{split}
\end{equation}
where
\begin{equation*}
\begin{split}
B(\delta)&=(\lambda^{p}-1)\left(\sum_{k=1}^{N}\alpha_{k}\big(h(P,u_{k})-\xi(P_{\delta})\cdot u_{k}\big)^{p}\right)\\
&\qquad+\alpha_{i_{0}}\lambda^{p}\bigg[\big(h(P,u_{i_{0}})-\xi(P_{\delta})\cdot
u_{i_{0}}-\delta\big)^{p}-\big(h(P,u_{i_{0}})-\xi(P_{\delta})\cdot
u_{i_{0}}\big)^{p}\bigg]\\
&=\left[(1-(c_{n}\delta^{n}+...+c_{2}\delta^{2}))^{-\frac{p}{n}}-1\right]\left(\sum_{k=1}^{N}\alpha_{k}(h(P,u_{k})-\xi(P_{\delta})\cdot u_{k})^{p}\right)\\
&\qquad+\alpha_{i_{0}}\lambda^{p}\bigg[\big(h(P,u_{i_{0}})-\xi(P_{\delta})\cdot
u_{i_{0}}-\delta\big)^{p}-\big(h(P,u_{i_{0}})-\xi(P_{\delta})\cdot u_{i_{0}}\big)^{p}\bigg].
\end{split}
\end{equation*}

From the facts that $d_{0}=d(P)>h(P,u_{i_{0}})-\xi(P_{\delta})\cdot u_{i_{0}}>h(P,u_{i_{0}})-\xi(P_{\delta})\cdot u_{i_{0}}-\delta>0$, $p<0$ and the fact that $f(t)=t^{p}$ is convex on $(0,\infty)$, we have
$$
\big(h(P,u_{i_{0}})-\xi(P_{\delta})\cdot u_{i_{0}}-\delta\big)^{p}-\big(h(P,u_{i_{0}})-\xi(P_{\delta})\cdot
u_{i_{0}}\big)^{p}>(d_{0}-\delta)^{p}-d_{0}^{p}>0.
$$
Hence,
\begin{equation}\tag{4.6}\label{Equation 4.6}
\begin{split}
B(\delta)&=(\lambda^{p}-1)\left(\sum_{k=1}^{N}\alpha_{k}\big(h(P,u_{k})-\xi(P_{\delta})\cdot u_{k}\big)^{p}\right)\\
&\qquad+\alpha_{i_{0}}\lambda^{p}\bigg[\big(h(P,u_{i_{0}})-\xi(P_{\delta})\cdot
u_{i_{0}}-\delta\big)^{p}-\big(h(P,u_{i_{0}})-\xi(P_{\delta})\cdot
u_{i_{0}}\big)^{p}\bigg]\\
&>\left[(1-(c_{n}\delta^{n}+...+c_{2}\delta^{2}))^{-\frac{p}{n}}-1\right]\left(\sum_{k=1}^{N}\alpha_{k}(h(P,u_{k})-\xi(P_{\delta})\cdot
u_{k})^{p}\right)\\
&\qquad+\alpha_{i_{0}}\lambda^{p}\big[(d_{0}-\delta)^{p}-d_{0}^{p}\big].
\end{split}
\end{equation}
On the other hand,
\begin{equation}\tag{4.7}\label{Equation 4.7}
\lim_{\delta\rightarrow0}\sum_{k=1}^{N}\alpha_{k}\big(h(P,u_{k})-\xi(P_{\delta})\cdot u_{k}\big)^{p}=\sum_{k=1}^{N}\alpha_{k}h(P,u_{k})^{p},
\end{equation}
\begin{equation}\tag{4.8}\label{Equation 4.8}
(d_{0}-\delta)^{p}-d_{0}^{p}>0,
\end{equation}
and
\begin{equation}\tag{4.9}\label{Equation 4.9}
\begin{split}
\lim_{\delta\rightarrow0}&\frac{(1-(c_{n}\delta^{n}+...+c_{2}\delta^{2}))^{-\frac{p}{n}}-1}{(d_{0}-\delta)^{p}-d_{0}^{p}}\\
&=\lim_{\delta\rightarrow0}\frac{(-\frac{p}{n})(1-(c_{n}\delta^{n}+...+c_{2}\delta^{2}))^{-\frac{p}{n}-1}(-nc_{n}\delta^{n-1}-...-2c_{2}\delta)}{p(d_{0}-\delta)^{p-1}(-1)}=0.
\end{split}
\end{equation}

From Equations (\ref{Equation 4.6}), (\ref{Equation 4.7}), (\ref{Equation 4.8}), (\ref{Equation 4.9}), and the fact that $p<0$, we have $B(\delta)>0$ for small enough $\delta>0$. From this and Equation (\ref{Equation 4.5}), there
exists a $\delta_{0}>0$ such that $P_{\delta_{0}}\in\mathcal{P}(u_{1},...,u_{N})$ and
$$
\Phi_{\lambda_{0}P_{\delta_{0}}}(\xi(\lambda_{0}P_{\delta_{0}}))>\Phi_{P}(\xi(P_{\delta_{0}}))\geq\Phi_{P}(\xi(P))=\Phi_{P}(o),
$$
where $\lambda_{0}=V(P_{\delta_{0}})^{-\frac{1}{n}}$. Let $P_{0}=\lambda_{0}P_{\delta_{0}}-\xi(\lambda_{0}P_{\delta_{0}})$, then
$P_{0}\in\mathcal{P}^{n}(u_{1},...,u_{N})$, $V(P_{0})=1$, $\xi(P_{0})=o$ and
\begin{equation}\tag{4.10}\label{Equation 4.10}
\Phi_{P_{0}}(o)<\Phi_{P}(o).
\end{equation}
This contradicts Equation (\ref{Equation 4.4}). Therefore, $P\in\mathcal{P}_{N}(u_{1},...,u_{N})$.
\end{proof}

Now we have prepared enough to prove the main theorem of this paper. We only need to prove the following:
\begin{theorem}
Suppose $p<0$, $\alpha_{1},...,\alpha_{N}>0$, and the unit vectors $u_{1},...,u_{N}$ are not concentrated on a hemisphere. If for any subspace $X$ with $1\leq\dim X\leq n-1$, $\{u_{1},...,u_{N}\}\cap X$ is always concentrated on a closed hemisphere of $S^{n-1}\cap X$,
then there exists a polytope $P_{0}\in\mathcal{P}_{N}(u_{1},...,u_{N})$ such that
$$
S_{p}(P_{0},\cdot)=\sum_{k=1}^{N}\alpha_{k}\delta_{u_{k}}(\cdot).
$$
\end{theorem}
\begin{proof}
Theorem 4.5 can be directly got by Lemma 3.4 and Lemma 4.4.
\end{proof}

\end{document}